\documentclass[a4paper]{article}

\usepackage{mathtools}
\usepackage{amsthm}
\usepackage{amssymb}
\usepackage{mathrsfs}

\DeclareMathOperator{\car}{char}
\DeclareMathOperator{\codim}{codim}
\DeclareMathOperator{\dvol}{dvol}
\DeclareMathOperator{\ind}{ind}
\DeclareMathOperator{\trdeg}{trdeg}
\DeclareMathOperator{\vol}{vol}

\DeclareMathOperator{\GL}{GL}
\DeclareMathOperator{\Spec}{Spec}

\DeclarePairedDelimiter{\abs}{\lvert}{\rvert}
\DeclarePairedDelimiter{\floor}{\lfloor}{\rfloor}
\DeclarePairedDelimiter{\norm}{\lVert}{\rVert}

\newtheorem{theorem}{Theorem}[section]
\newtheorem{corollary}[theorem]{Corollary}
\newtheorem{hypothesis}[theorem]{Hypothesis}
\newtheorem{lemma}[theorem]{Lemma}
\newtheorem{proposition}[theorem]{Proposition}

\theoremstyle{remark}
\newtheorem{remark}[theorem]{Remark}

\newcommand{\bC}{{\mathbb C}}
\newcommand{\bN}{{\mathbb N}}
\newcommand{\bP}{{\mathbb P}}
\newcommand{\bQ}{{\mathbb Q}}
\newcommand{\bR}{{\mathbb R}}
\newcommand{\bZ}{{\mathbb Z}}

\newcommand{\cA}{{\mathcal A}}
\newcommand{\cB}{{\mathcal B}}
\newcommand{\cF}{{\mathcal F}}
\newcommand{\cH}{{\mathcal H}}
\newcommand{\cM}{{\mathcal M}}
\newcommand{\cN}{{\mathcal N}}
\newcommand{\cO}{{\mathcal O}}
\newcommand{\cP}{{\mathcal P}}
\newcommand{\cQ}{{\mathcal Q}}
\newcommand{\cR}{{\mathcal R}}
\newcommand{\cX}{{\mathcal X}}
\newcommand{\cY}{{\mathcal Y}}

\newcommand{\sD}{{\mathscr D}}

\newcommand{\Kbar}{{\overline K}}

\title{A generalized Vojta--R\'emond inequality over adelic curves}
\author{Ningjun Jiang\thanks{The author would like to thank Professor Huayi Chen for introducing him to Arakelov theory of adelic curves.}}

\date{November 2026}

\begin{document}

\maketitle

\begin{abstract}
We extend Dill's generalized Vojta--R\'emond inequality to proper adelic curves of characteristic zero, under an explicit Siegel hypothesis. The geometric case is closely related to recent work on Mordell--Lang plus geometric Bogomolov for semi-abelian varieties. Our formulation retains Dill's generality for products of distinct projective varieties. The appendix establishes the Siegel input for arithmetic function fields with nef and big polarizations on the model. 
\end{abstract}

\tableofcontents

\section{Introduction}\label{sec:introduction}

R\'emond's generalized Vojta inequality grew out of height inequalities used in proofs of the Mordell--Lang conjecture. It was used in the proof of Poonen's Mordell--Lang plus Bogomolov conjecture for semi-abelian varieties over number fields \cite[\S1.1]{Remond2005}. In its general form, it compares two heights on a modification of a power of a projective variety, under a lower bound for certain geometric intersection numbers \cite[Th\'eor\`eme~1.2]{Remond2005}. Dill extended this comparison to products of possibly different projective varieties, with individual height thresholds for the factors \cite[Theorem~1.1]{Dill2020}. This generalization was motivated by a relative Mordell-–Lang problem involving isogeny orbits in an abelian scheme, where the factors may lie in different fibers \cite[\S1]{Dill2020}.

Over geometric function fields of characteristic zero, Long Liu proves a Vojta inequality for semi-abelian varieties, excluding tuples of sufficiently large height satisfying both angular and spacing conditions \cite[Theorem~4.1]{Liu2026}. For this he proves a geometric version of \cite[Th\'eor\`eme~1.2]{Remond2005} over the function field of a curve. This inequality is an input to his proof of Mordell--Lang plus geometric Bogomolov. Similar results were established independently in the author's thesis \cite{Jiang2026}. The two approaches differ in their treatment of higher-dimensional bases: Liu reduces to the curve case using a generic complete-intersection curve \cite[\S4.8.2 and Lemma~4.20]{Liu2026}, while the thesis constructs auxiliary sections by applying a Siegel lemma directly on higher-dimensional polarized models.

Following this direct approach, we extend the general height comparison to proper adelic curves in the sense of Chen and Moriwaki \cite{ChenMoriwaki2020}, unifying the arithmetic and geometric scenarios above. We retain Dill's products of possibly different projective varieties and individual height thresholds. The construction of auxiliary sections requires a bound for small solutions of linear systems, which we formulate as a Siegel hypothesis on the adelic curve. This gives the comparison over any proper adelic curve satisfying that hypothesis. We work throughout in characteristic zero; the Frobenius example in \cite[\S5.3]{Liu2026} explains why an unrestricted extension of this approach to positive characteristic is not expected.

Let
\[ S=(K,(\Omega,\cA,\nu),\phi), \qquad \phi(v)=\abs\cdot_v,\]
be a proper adelic curve with $\car K=0$, and fix an algebraic closure $\Kbar$ of $K$. Thus, for every $b\in K^\times$, the function $v\mapsto\log\abs b_v$ is measurable and integrable, and
\begin{equation}\label{eq:product-formula}
  \int_\Omega\log\abs b_v\,d\nu=0.
\end{equation}
We assume throughout that $\abs\cdot_v$ is the usual absolute value on $K_v=\bR$ or $\bC$ for $v\in\Omega_\infty$, but see Remark~\ref{rem:normalize-adelic-curve}. We use the algebraic extensions $S_F=S\otimes_K F$ in \cite[\S\S3.3--3.4]{ChenMoriwaki2020}. In particular, the measures on these extensions incorporate the normalization by the field degree. For $z=[z_0:\dotsb:z_N]\in\bP^N(\Kbar)$, choose a finite extension $F/K$ containing its coordinates and set
\begin{equation}\label{eq:point-height}
h(z)=h_S(z)=\int_{\Omega_F}\log\max_{0\leqslant j\leqslant N}\abs {z_j}_v\,d\nu_F.
\end{equation}
This is independent of the coordinates chosen and of $F$, by \eqref{eq:product-formula} and compatibility of algebraic extensions \cite[Theorem~3.4.12 and Definition~3.5.1]{ChenMoriwaki2020}. In the following we suppress the fixed adelic curve $S$ from the height notation. For a nonzero polynomial, $h$ denotes the height of its coefficient vector. For a family of polynomials, it denotes the height of their joint coefficient vector.

With this height convention, the input needed for the auxiliary-section construction takes the following form. The section is represented by a nonzero vector in the kernel of a linear system, and its height must be controlled in terms of the heights of the equations. In the systems arising in the proof, the ratio of the number of unknowns to the dimension of the kernel stays bounded. We therefore require a Siegel estimate uniform under such a bound, over each finite extension of $K$.

\begin{hypothesis}[Siegel estimate]\label{hyp:siegel}
There exists $\sigma_S\geqslant 1$ such that, for every finite extension $F/K$ and every $R>1$, there are constants $\gamma_S(F,R),c_S(F,R)\geqslant 0$ with the following property: for every matrix $A\in M_{b\times n}(F)$ of rank $b$, with $1\leqslant b<n$ and $n/(n-b)\leqslant R$, there is $0\neq z\in F^n$ such that $Az=0$ and
\begin{equation}\label{eq:siegel-hypothesis}
h(z)\leqslant\frac{n}{n-b} \bigl(\sigma_S\max_{1\leqslant j\leqslant b}h(A_j)+\gamma_S(F,R)\log n+c_S(F,R)\bigr),
\end{equation}
where $A_j$ is the $j$-th row of $A$.
\end{hypothesis}

The hypothesis is verified for number fields in \cite[p.~7]{Bombieri1979} with $(\sigma_S,\gamma_S)=(1,1)$ and for geometric function fields in \cite{Jiang2026} with $(\sigma_S,\gamma_S)=(\trdeg_kK,0)$. The appendix verifies the hypothesis for arithmetic function fields with nef and big polarizations on the model \cite{Moriwaki2000}. For a general proper adelic curve, it remains an additional assumption. 

We now introduce the geometric data for the height comparison as in \cite{Dill2020}. Let $m\geqslant 2$, and fix closed embeddings $X_i\subset\bP^{N_i}_{\Kbar}$ of integral projective varieties of positive dimension, with homogeneous coordinates
\[ W^{(i)}=(W^{(i)}_0,\dotsc,W^{(i)}_{N_i}) \qquad(1\leqslant i\leqslant m).\]
All point heights, degrees, and subvariety heights on $X_i$ refer to this fixed embedding; in particular, we write $h(x_i)$ for the height of $x_i\in X_i(\Kbar)\subset\bP^{N_i}(\Kbar)$. Write $X=X_1\times\dotsb\times X_m$ and let $p_i\colon X\to X_i$ be projections. For $a=(a_1,\dotsc,a_m)\in(\bN_+)^m$, put
\[ \cN_a=\bigotimes_{i=1}^m p_i^*\cO_{X_i}(a_i).\]
Fix a nonempty open subset $U^0\subset X$. We consider triples $(a,\cX,\cM)$, where $\cX$ is an integral projective variety containing $U^0$ as an open subset, equipped with a proper morphism $\pi\colon \cX\to X$ restricting to the identity on $U^0$, and $\cM$ is a nef line bundle on $\cX$. We also write $\cN_a$ for $\pi^*\cN_a$ and $W^{(i)}$ for pulled-back coordinates.

Fix integers $t_1,t_2,M\geqslant 1$ and real numbers $\delta_1,\dotsc,\delta_m\geqslant 1$, independently of the triple. We assume that there exists a very ample line bundle $\cP$ on $\cX$ and injections
\begin{equation}\label{eq:injections}
\cP\hookrightarrow\cN_a^{\otimes t_1},\qquad \cP\otimes\cM^{-1}\hookrightarrow\cN_a^{\otimes t_2},
\end{equation}
both of which restrict to isomorphisms on $U^0$, with the following properties: there is a system of homogeneous coordinates $\Xi$ for $\cP$ whose images under the first injection are unit monomials of multidegree $t_1a$ in $W^{(i)}$; there is also a family of sections $Z=(Z_1,\dotsc,Z_M)$ generating $\cP\otimes\cM^{-1}$ on $\cX$, whose images under the second injection are multihomogeneous polynomials $P_1,\dotsc,P_M$ of multidegree $t_2a$ satisfying
\begin{equation}\label{eq:coefficient-bound}
 h(P_1,\dotsc,P_M)\leqslant\sum_{i=1}^m a_i\delta_i.
\end{equation}
All these data descend to a finite extension of $K$, which may be enlarged without changing the heights.

For $x=(x_1,\dotsc,x_m)\in U^0(\Kbar)$, define
\[\begin{split}
  h_{\cM}(x)&=h(\Xi(x))-h(Z(x)),\label{eq:M-height}\\
  h_{\cN_a}(x)&=\sum_{i=1}^m a_i h(x_i).
\end{split}\]

Let $\theta\geqslant1$ and $\omega\geqslant-1$ be integers, and set
\[ u_0=\sum_{i=1}^m\dim X_i, \qquad \omega'=3+\omega.\]
Put $\beta_S=\nu(\Omega_\infty)$ and $\delta_i^*=\max(\beta_S,\delta_i)$, and let $\sigma_S\geqslant 1$ be the constant in Hypothesis~\ref{hyp:siegel}. We keep Dill's constant:
\begin{equation}\label{eq:Lambda}
 \Lambda=\theta(2t_1u_0)^{u_0}\bigl(\max_{1\leqslant i\leqslant m}N_i+1\bigr)\prod_{i=1}^m\deg(X_i).
\end{equation}
For $0\leqslant u\leqslant u_0$ define
\[\psi(u)=\prod_{j=u+1}^{u_0}(\omega'j+1),\qquad c_1=c_2=\Lambda^{\psi(0)},\]
and put $B=\sigma_SMt_2$ and
\[c_3^{(i)}=\Lambda^{2\psi(0)}B^{u_0}\bigl(h(X_i)+\delta_i^*\bigr)\qquad(1\leqslant i\leqslant m),\]
where $h(X_i)$ is the height of $X_i$ inside $\bP^{N_i}$ in the sense of \cite[Definition~4.4.3]{ChenMoriwaki2021}. We specify this height convention below.

\begin{theorem}\label{thm:vojta-remond}
Assume Hypothesis~\ref{hyp:siegel}. Let $(a,\cX,\cM)$ satisfy the conditions above, and let $x=(x_1,\dotsc,x_m)\in U^0(\Kbar)$. Suppose that, for every product
\[Y=Y_1\times\dotsb\times Y_m\]
with $Y_i\subset X_i$ a closed integral subvariety containing $x_i$, one has
\begin{equation}\label{eq:intersection-bound}
 \bigl(\cM^{\dim Y}\cdot\cY\bigr)\geqslant\theta^{-1}\prod_{i=1}^m\deg(Y_i)^{-\omega}a_i^{\dim Y_i},
\end{equation}
where $\cY$ is the closure of $\pi^{-1}(Y\cap U^0)$ in $\cX$. If
\begin{equation}\label{eq:height-hypotheses}
 c_2a_{i+1}\leqslant a_i\quad(1\leqslant i<m),\qquad c_3^{(i)}\leqslant h(x_i)\quad(1\leqslant i\leqslant m),
\end{equation}
then
\begin{equation}\label{eq:main-inequality}
 h_{\cN_a}(x)\leqslant c_1h_{\cM}(x).
\end{equation}
\end{theorem}

For the usual adelic structure on a number field, we have $\beta_S=\sigma_S=1$, and the constants and conclusion agree with \cite[Theorem~1.1]{Dill2020}. For a geometric function field equipped with an ample polarization on a projective model, the theorem gives a height comparison directly over that field. Over the function field of a curve, taking all $X_i$ equal, all $\delta_i^*=\delta_i$ equal, and $\omega=0$ gives a comparison of the form stated in \cite[Proposition~4.2]{Liu2026}.

We now fix the convention for heights of subvarieties. For an integral closed subvariety $V\subset\bP^N_{\Kbar}$ of dimension $d$, let $\deg(V)=(\cO(1)^d\cdot V)$, and let $h(V)$ be the arithmetic intersection number \cite[Definition~4.4.3]{ChenMoriwaki2021}
\begin{equation}\label{eq:variety-height}
 h(V)=\bigl(\overline{\cO(1)}\mid V\bigr)^{d+1}_{S_F},
\end{equation}
computed over a finite extension $F/K$ where everything is defined. Here $\overline{\cO(1)}$ is the adelic line bundle equipped with the orthogonal quotient metric family induced by the standard norm family on $F^{N+1}$ \cite[Examples~4.1.2 and 4.1.7]{ChenMoriwaki2021}. Compatibility with extension of scalars follows from \cite[Theorem~4.3.6]{ChenMoriwaki2021}. 

For the same choice of Chow form, the local terms in \cite[Theorem~3.8.7]{ChenMoriwaki2021} agree with those defining the height in \cite[\S2.1, pp.~96--97]{Remond2001Multiprojective}: the logarithm of the maximum coefficient norm at nonarchimedean places, and the spherical average with the same harmonic correction at archimedean places. The product formula removes the choice of scalar for the Chow form. Thus this convention recovers exactly the subvariety height used by Dill in the number field case. 

For a point $z$, the linear-form identity in \cite[\S2.1, p.~96]{Remond2001Multiprojective} shows that $h(\{z\})$ uses Euclidean norms at archimedean places. Since the point height in \eqref{eq:point-height} uses maximum norms, we have $h(z)\leqslant h(\{z\})$.

The total archimedean mass $\beta_S=\nu(\Omega_\infty)$ of $S$ is finite by \cite[Proposition~3.1.2]{ChenMoriwaki2020}. It records the archimedean contribution of integer constants and the number of summands when local height estimates are integrated; it equals $1$ for number fields and $0$ for geometric function fields.

\begin{remark}\label{rem:normalize-adelic-curve}
Our assumption on archimedean absolute values entails no loss of generality. As explained at the beginning of \cite[\S4.1]{ChenMoriwaki2020}, if originally $\abs\cdot_v=\abs\cdot_{\mathrm{usual}}^{\kappa(v)}$, replace it by the usual absolute value and replace $d\nu$ by $\kappa(v)d\nu$ on $\Omega_\infty$, leaving nonarchimedean data unchanged. This preserves the product formula, commutes with finite algebraic extension, and leaves all point and coefficient heights unchanged. For an unnormalized adelic curve, subvariety heights are understood via this normalization and the archimedean factor is $\int_{\Omega_\infty}\kappa(v)\,d\nu$. With these conventions, Hypothesis~\ref{hyp:siegel} and Theorem~\ref{thm:vojta-remond} retain the same constants for the original adelic curve.
\end{remark}

The proof is organized as in \cite{Dill2020}. It follows R\'emond's reduction to a product of minimal dimension, together with auxiliary-section and product arguments \cite{Remond2001Product,Remond2005}. R\'emond's calculations with local absolute values and local heights, including the multiprojective estimates in \cite{Remond2001Multiprojective}, are used place by place. For the passage to global heights we use the adelic formalism of \cite{ChenMoriwaki2020,ChenMoriwaki2021}. The existence of a sufficiently small nonzero solution of a linear system needed for constructing the auxiliary section is a separate input, i.e., Hypothesis~\ref{hyp:siegel}.

\paragraph{AI disclosure}

During the preparation of this work, large language models (LLMs) were used to assist in locating and explaining references and to support the writing process.

\section{Reduction to a minimal subproduct}\label{sec:reduction}

Fix a triple $(a,\cX,\cM)$ and a point $x$ as in Theorem~\ref{thm:vojta-remond}, satisfying \eqref{eq:height-hypotheses}. We follow \cite[\S2]{Dill2020}, keeping track of the archimedean factor $\beta_S$. The value of $\Lambda$ in \eqref{eq:Lambda}, or any larger value, suffices for all the estimates below.

We first record the form of arithmetic B\'ezout that will be used. The intersection formulas of \cite[Chapter~4]{ChenMoriwaki2021} are stated under a countability or discreteness assumption. For coordinate metrics considered here, we apply them after choosing a finite extension $F/K$ containing the data and descending their finitely many defining coefficients to a subfield $F_0\subset F$ which is finitely generated over $\bQ$. Restrict $S_F$ to $F_0$ as in \cite[\S3.2.8]{ChenMoriwaki2020} and denote the restriction by $S_0$. This is a proper adelic curve with countable underlying field and the same parameter space and measure. The local coordinate heights are unchanged by this restriction: for each $v\in\Omega_F$, the inclusion $F_0\subset F$ induces an isometric embedding $(F_0)_v\hookrightarrow F_v$ of completions, and these heights are compatible with extension of scalars; see \cite[Proposition~3.8.9 and \S3.9]{ChenMoriwaki2021} for the nonarchimedean case, the archimedean case being immediate. Thus the identities so obtained apply to the original data. The archimedean mass $\beta_S=\beta_{S_F}=\beta_{S_0}$ remains unchanged.

\begin{lemma}[Heights of hypersurface sections]\label{lem:bezout}
Let $T\subset\bP^N_{\Kbar}$ be an integral projective variety of positive dimension, and let $F$ be a homogeneous polynomial of degree $e\geqslant 1$ whose restriction to $T$ is nonzero. If $T'$ is an irreducible component of $T\cap V(F)$, with its reduced structure, then
\[\begin{aligned}
\deg(T')&\leqslant e\deg(T),\\
h(T')&\leqslant e h(T)+\deg(T)\bigl(h(F)+\beta_S\sqrt N\bigr).
\end{aligned}\]
\end{lemma}

\begin{proof}
Write $Z=T\cdot V(F)$ for the intersection cycle, including multiplicities, and $r=\dim T$. With $\overline L=\overline{\cO(1)}\lvert_T$, \cite[Proposition~4.4.4(3)]{ChenMoriwaki2021} gives, over a suitable finite extension $K'/K$,
\begin{equation}\label{eq:section-height-identity}
h(Z)=e h(T)+ \int_{\Omega_{K'}}\Bigl(\int_{T_v^{\mathrm{an}}}\log\abs F_{\overline L^{\otimes e},v}\,d\mu_v\Bigr)\,d\nu_{K'},
\end{equation}
where $\mu_v$ is the positive local intersection measure associated with $r$ copies of $\overline L_v$, of mass $\deg(T)$. At a nonarchimedean place the section norm is at most the maximum absolute value of its coefficients. At an archimedean place, writing $F=\sum_{\abs\alpha=e}c_\alpha X^\alpha$, Cauchy--Schwarz gives
\[
\sup\abs F_{\overline L^{\otimes e},v} \leqslant\bigl(\sum_{\abs\alpha=e}\frac{\alpha!}{e!}\abs{c_\alpha}_v^2\bigr)^{1/2}\leqslant\mathrm e^{\sqrt N}\max_\alpha\abs{c_\alpha}_v.
\]
The last inequality is \cite[Lemme~5.2]{Remond2001Product}. Then the integration in \eqref{eq:section-height-identity} proves the height bound for $Z$. Heights of subvarieties are nonnegative: this follows by induction on dimension from \cite[Proposition~4.4.4(3)]{ChenMoriwaki2021}, using a coordinate section of norm at most $1$, and the nonnegativity of point heights. Hence each component has height at most $h(Z)$ by \cite[Proposition~4.4.4(2)]{ChenMoriwaki2021}. Finally $\deg(Z)=e\deg(T)$, which proves the degree bound as well.
\end{proof}

Write
\[\begin{aligned}
 D_i&=\deg(X_i),\qquad &A&=\max_i(N_i+1),\\
 \sD&=\prod_iD_i,\qquad &H_0&=\sum_i a_i\bigl(h(X_i)+\delta_i^*\bigr)>0.
\end{aligned}\]
Choose a product $Y=Y_1\times\dotsb\times Y_m$ of closed integral subvarieties of the $X_i$, of minimal total dimension $u=\sum_i u_i$, where $u_i=\dim Y_i$, subject to
\begin{align}
  x_i&\in Y_i(\Kbar) &(1\leqslant i\leqslant m),\label{eq:minimal-points}\\
  d_i&\leqslant D_i\Lambda^{\psi(u)-1} &(1\leqslant i\leqslant m),\label{eq:minimal-degrees}\\
  \prod_i d_i&\leqslant\sD\Lambda^{\psi(u)-1},\label{eq:minimal-product}\\
  H(Y)\coloneqq\sum_i a_i(h_i+\delta_i^*)&\leqslant\frac12\Lambda^{2\psi(u)}B^{u_0-u}H_0,\label{eq:minimal-heights}
\end{align}
with $d_i=\deg(Y_i)$ and $h_i=h(Y_i)$ computed in the embedding $Y_i\hookrightarrow X_i\hookrightarrow\bP^{N_i}$. Such a product exists: take $Y=X$ and use $\psi(u_0)=1$ and $\Lambda\geqslant 2$. Moreover we have $u>0$. Otherwise $Y_i=\{x_i\}$ for every $i$, and \eqref{eq:height-hypotheses} and \eqref{eq:minimal-heights} give
\[ 0<\sum_i a_i c_3^{(i)}\leqslant h_{\cN_a}(x)\leqslant\sum_i a_i h_i\leqslant\frac12\sum_i a_i c_3^{(i)},\]
a contradiction.

For this product define
\begin{equation}\label{eq:Lambda-h}
 \Lambda_h=\sum_i a_i\bigl(h_i+\delta_i^*+\beta_S d_i(u_i+1)\log(2d_i(N_i+1))\bigr).
\end{equation}
We claim that
\begin{equation}\label{eq:Lambda-h-bound}
 \Lambda_h<\Lambda^{2\psi(u)}B^{u_0-u}H_0 =\Lambda^{2\psi(u)-2\psi(0)}B^{-u}\sum_i a_i c_3^{(i)}.
\end{equation}
Indeed, $H_0\geqslant\beta_S\sum_i a_i$ and $\log t<t$ give
\[\begin{aligned}
 \Lambda_h-H(Y) &\leqslant H_0\sum_i d_i(u_i+1)\log\bigl(2d_i(N_i+1)\bigr)\\
 &\leqslant 2(u+m)A\sD^2\Lambda^{2\psi(u)-2}H_0 <\frac12\Lambda^{2\psi(u)}H_0.
\end{aligned}\]
The last inequality follows from $u+m\leqslant 2u_0$ and
\[\Lambda^2\geqslant (2u_0)^{2u_0}A^2\sD^2>8u_0A\sD^2.\]
Together with \eqref{eq:minimal-heights}, this proves
\eqref{eq:Lambda-h-bound}. 

We next choose adapted coordinates, in the sense of \cite[\S2.1]{Remond2005}. The algebraic construction of \cite[Proposition~2.2]{Remond2005}, valid over any algebraically closed field of characteristic zero, gives
\[ V_j^{(i)}=\sum_{k=0}^{N_i}M_{jk}^{(i)}W_k^{(i)},\]
with
\[M^{(i)}\in\GL_{N_i+1}(\bQ)\cap M_{N_i+1}(\bZ),\quad\abs{M_{jk}^{(i)}}\leqslant\max\{1,d_i/2\}.\]
In particular, the first $u_i+1$ coordinates define a finite surjective morphism $\rho_i:Y_i\to\bP^{u_i}$ of degree $d_i$. If $u_i<N_i$, they may be chosen so that $V_{u_i+1}^{(i)}/V_0^{(i)}$ generates the function field of $Y_i$ over that of $\bP^{u_i}$. If $Y_i=\bP^{N_i}$, we take $M^{(i)}$ to be the identity.

\begin{proposition}\label{prop:no-small-equation}
There is no pair $(l,U)$, where $1\leqslant l\leqslant m$ and $U$ is a nonzero homogeneous $\Kbar$-polynomial in $(V_0^{(l)},\dotsc,V_{u_l}^{(l)})$ such that
\begin{align}
U(x_l)&=0,\label{eq:small-equation-zero}\\
\deg U&\leqslant\Lambda^{\omega'u\psi(u)},\label{eq:small-equation-degree}\\
a_l h(U)&\leqslant\Lambda^{2\psi(u-1)-2\psi(u)}\frac{B}{4d_l}\Lambda_h.\label{eq:small-equation-height}
\end{align}
\end{proposition}

\begin{proof}
Suppose such a pair exists and put $e=\deg U\geqslant1$ and $E=\Lambda^{\omega'u\psi(u)}$. The finite surjective map $\rho_l$ ensures that $U(V^{(l)})$ does not vanish identically on $Y_l$. Choose an irreducible component $Y_l'$ of its zero locus containing $x_l$. Set $Y_i'=Y_i$ for $i\neq l$ and put $Y'=\prod_iY_i'$. We have $u_l\geqslant 1$ and $\dim Y'=u-1$.

Express $U(V^{(l)})$ in the original coordinates as $U'(W^{(l)})$. The integer entries of $M^{(l)}$ have absolute value at most $1$ at nonarchimedean places. The same coefficient calculation as in \cite[Proposition~3.1]{Remond2005} and \cite[(2.2)]{Dill2020}, integrated with the archimedean weights, yields
\begin{equation}\label{eq:coordinate-height-change}
\begin{split}
 h(U')&\leqslant h(U)+\beta_S\biggl(e\log\bigl((N_l+1)\max\{1,d_l/2\}\bigr)+\log\binom{e+u_l}{e}\biggr)\\
 &\leqslant h(U)+\beta_S e\log\bigl(d_l(N_l+1)(u_l+1)\bigr).
\end{split}
\end{equation}
Here we use
\[\binom{e+u_l}{e}=\prod_{i=1}^e\bigl(1+\frac{u_l}{i}\bigr)\leqslant(u_l+1)^e.\]
Lemma~\ref{lem:bezout} therefore gives
\[ \deg(Y_l')\leqslant e d_l,\qquad h(Y_l')\leqslant e h_l+d_l\bigl(h(U')+\beta_S\sqrt{N_l}\bigr).\]
Since
\begin{equation}\label{eq:psi-recurrence}
 \psi(u-1)-\psi(u)=\omega'u\psi(u),
\end{equation}
the degree conditions \eqref{eq:minimal-degrees} and
\eqref{eq:minimal-product} hold for $Y'$.

For the height condition, $e\geqslant 1$ and the nonnegativity of $h_i$ give
\[\begin{split}
 H(Y')\leqslant{} &eH(Y)+a_l d_l\Bigl(h(U)+\beta_S\Bigl(e\log\bigl(d_l(N_l+1)(u_l+1)\bigr)+\sqrt{N_l}\Bigr)\Bigr)\\
 \leqslant{} &a_l d_l h(U)+3E\Lambda_h.
\end{split}\]
For the last step, bound the first term involving $\beta_S$ using $e\leqslant E$ and \[\log(d_l(N_l+1)(u_l+1))\leqslant(u_l+1)\log(2d_l(N_l+1)),\]
and the second using $\beta_S a_l d_l\leqslant\Lambda_h$ and $\sqrt{N_l}\leqslant E$. Now $E\geqslant\Lambda^2\geqslant 12$ and $B\geqslant 1$, so $3E\leqslant BE^2/4$. By \eqref{eq:small-equation-height} and \eqref{eq:Lambda-h-bound},
\[ H(Y')\leqslant\frac12B E^2\Lambda_h\leqslant\frac12\Lambda^{2\psi(u-1)}B^{u_0-(u-1)}H_0.\]
Thus $Y'$ satisfies all four minimality conditions and has smaller dimension, a contradiction.
\end{proof}

We finish this section with consequences of the proposition for local coordinates as in \cite[Corollary~2.2]{Dill2020}. The construction of \cite[Lemmes~2.1 and~2.3]{Remond2005} supplies, for $0\leqslant j\leqslant N_i$, homogeneous polynomials $P_j^{(i)},Q_j^{(i)}\in\Kbar[T_0,\dotsc,T_{u_i},T_{u_i+1}]$ of degree $d_i$, monic in the last variable and each a power of an irreducible polynomial, such that
\[ P_j^{(i)}(V_0^{(i)},\dotsc,V_{u_i}^{(i)},V_j^{(i)})=Q_j^{(i)}(V_0^{(i)},\dotsc,V_{u_i}^{(i)},W_j^{(i)})=0 \quad\text{on }Y_i.\]
Let $b_{i,j}$ be the exponent of the irreducible factor of $Q_j^{(i)}$, and let $B_i$ be the joint coefficient vector of all these polynomials. The local calculation in the proof of that lemma gives
\begin{equation}\label{eq:relation-coefficients}
 h(B_i)\leqslant h_i+\beta_S d_i(u_i+1)\log\bigl(d_i(N_i+1)\bigr).
\end{equation}
More explicitly, let $f_i$ be the Chow form in the original coordinates, normalized by $f_i(V_0^{(i)},\dotsc,V_{u_i}^{(i)})=1$. Write $M_v(f_i)$ for the corrected spherical measure of
\cite[\S2.1]{Remond2001Multiprojective} at an archimedean place, and for the maximum coefficient norm otherwise. By \cite[Theorem~3.8.7 and Remark~4.2.13]{ChenMoriwaki2021},
\[\int\log M_v(f_i)\,d\nu=h_i,\]
with the integral taken over a field of definition. The local estimate in \cite[proof of Lemme~2.3]{Remond2005}, using \cite[proof of Lemme~3.3]{Remond2001Multiprojective}, gives at an archimedean place
\[ \log\max_{b\in B_i}\abs b_v \leqslant \log M_v(f_i)+d_i(u_i+1)\log\bigl(d_i(N_i+1)\bigr).\]
At a nonarchimedean place the same bound holds without the additive term. Integration gives exactly \eqref{eq:relation-coefficients}.

\begin{corollary}\label{cor:local-coordinates}
For every $i$, the morphism $\rho_i$ is \'etale at $x_i$, and $V_0^{(i)}(x_i)\neq 0$. Moreover, if $W_j^{(i)}$ does not vanish identically on $Y_i$, then $W_j^{(i)}(x_i)\neq 0$ and
\[ \frac{\partial^{b_{i,j}}Q_j^{(i)}}{\partial T_{u_i+1}^{b_{i,j}}}\biggl(1,\frac{V_1^{(i)}}{V_0^{(i)}},\dotsc,\frac{V_{u_i}^{(i)}}{V_0^{(i)}},\frac{W_j^{(i)}}{V_0^{(i)}}\biggr)(x_i)\neq 0.\]
\end{corollary}

\begin{proof}
We use the discriminant and resultant argument of \cite[Corollaire~3.2]{Remond2005} with the correction in
\cite[Corollary~2.2]{Dill2020}.
The polynomial $T_0$ detects the vanishing of $V_0^{(i)}$. If $u_i<N_i$, the primitive-element property makes $P_{u_i+1}^{(i)}$ irreducible; its nonzero discriminant with respect to the last variable detects the possible failure of \'etaleness. In the case $u_i=N_i$, $Y_i=\bP^{N_i}$ and the \'etale assertion is automatic. The resultants of $Q_j^{(i)}$ with its $b_{i,j}$-th derivative with respect to the last variable are nonzero in characteristic zero and detect vanishing of the displayed derivative. Finally, when $W_j^{(i)}$ is not identically zero, $Q_j^{(i)}(T_0,\dotsc,T_{u_i},0)$ is a nonzero polynomial and detects $W_j^{(i)}(x_i)=0$.

If any assertion failed, these constructions would produce a nonzero homogeneous $U$ vanishing at $\rho_i(x_i)$, with
\[ \deg U\leqslant 2d_i^2,\qquad h(U)\leqslant 6\beta_S N_i d_i^3+2d_i h(B_i).\]
These are precisely the determinant estimates in \cite[proof of Corollaire~3.2]{Remond2005}.\footnote{In R\'emond's notation, the printed binomial coefficient $\binom{D_i}{u_i}$ should read $\binom{D_i+u_i}{u_i}$, the number of degree-$D_i$ monomials in $u_i+1$ variables. This correction leaves the final bound unchanged.} The combinatorial error terms occur only at archimedean places and hence acquire the factor $\beta_S$.

Finally, the bounds fit Proposition~\ref{prop:no-small-equation} with the present constants. Set $q=\psi(u)$ and $E=\Lambda^{\omega'u q}\geqslant\Lambda^{2q}$. By \eqref{eq:minimal-degrees}, $d_i\leqslant\sD\Lambda^{q-1}$, whence $2d_i^2\leqslant E$. Also \eqref{eq:Lambda-h} and \eqref{eq:relation-coefficients} give $a_i h(B_i)\leqslant\Lambda_h$ and $\beta_S a_i\leqslant\Lambda_h$. Consequently
\[ 4d_i a_i h(U) \leqslant(24N_i d_i^4+8d_i^2)\Lambda_h\leqslant 32A\sD^4\Lambda^{4q-4}\Lambda_h\leqslant E^2B\Lambda_h,\]
where the last inequality uses $\Lambda^4\geqslant 32A\sD^4$. Together with \eqref{eq:psi-recurrence}, this contradicts Proposition~\ref{prop:no-small-equation}.
\end{proof}

\section{Constructing a section of small height}\label{sec:small-section}

We retain the minimal product $Y$ from Section~\ref{sec:reduction} and its strict transform $\cY\subset\cX$. Following \cite[\S4]{Remond2005} and \cite[\S3]{Dill2020}, we first obtain a lower bound for the dimension of a space of auxiliary sections, then describe this space as the kernel of a linear system of controlled height. Hypothesis~\ref{hyp:siegel}, stated in the introduction, supplies a nonzero section of small height from these estimates.

Fix a finite extension $F/K$ over which $\cY$, the line bundles, the chosen coordinates, and the polynomials in section~\ref{sec:reduction} are defined. Any further enlargement needed below is made once and is independent of $d$. Set
\begin{equation}\label{eq:epsilon}
 \varepsilon=\frac{1}{2u\theta(t_1m)^u}\prod_{i=1}^m d_i^{-1-\omega}, \qquad M'=M\prod_{i=1}^m(N_i+1).
\end{equation}
For integers $d>0$ with $d\varepsilon\in\bZ$, put
\[ \cQ_d=\cM^{\otimes d}\otimes\cN_a^{\otimes(-d\varepsilon)}, \qquad \cR_d=\cM^{\otimes(-d)}\otimes\cP^{\otimes d}\otimes\cN_a^{\otimes d\varepsilon},
\]
and let $Z'_d$ be the family
\[
 \Bigl(Z_j^{\otimes d}\otimes\bigotimes_{i=1}^m(W_{k_i}^{(i)})^{\otimes d\varepsilon a_i}\Bigr)_{\substack{1\leqslant j\leqslant M\\0\leqslant k_i\leqslant N_i}}.
\]
It has $M'$ entries and generates $\cR_d$. Restrict these bundles and sections to $\cY$.

By \cite[Proposition~4.1]{Remond2005}, with the nefness assumption as in \cite[proof of Proposition~3.1]{Dill2020}, for sufficiently large divisible $d$ one has
\begin{equation}\label{eq:section-dimension}
 k_d\coloneqq h^0(\cY,\cQ_d)\geqslant\frac{d^u}{4\theta u!}\prod_{i=1}^m d_i^{-\omega}a_i^{u_i}+o(d^u).
\end{equation}
This is a geometric estimate: the argument uses \eqref{eq:intersection-bound}, nefness of $\cM$, and
\begin{equation}\label{eq:intersection-upper-bound}
(\cN_a^u\cdot\cY)=\frac{u!}{u_1!\dotsm u_m!}\prod_i d_i a_i^{u_i} \leqslant m^u\prod_i d_i a_i^{u_i}.
\end{equation}
The required injection $\cM\hookrightarrow\cN_a^{\otimes t_1}$ on $\cY$ is obtained by choosing a member of $Z$ nonzero at $x$ and using the first injection in \eqref{eq:injections}. Choose $d_0$ divisible by the denominator of $\varepsilon$ and large enough so that, for all $d\in d_0\bN$, a monomial subspace $E_d\subset F[\Xi]_d$ maps isomorphically onto $H^0(\cY,\cP^{\otimes d})$. Such subspaces exist by Serre vanishing for the embedding given by $\Xi$. All subsequent asymptotics refer to $d\to\infty$ with other data fixed. 

\begin{proposition}\label{prop:small-section}
Assume Hypothesis~\ref{hyp:siegel}. For every sufficiently large $d\in d_0\bN$, there exists
$0\neq s_d\in H^0(\cY,\cQ_d)$ such that
\begin{equation}\label{eq:small-section-height}
 h(s_d)\leqslant\frac{2\sigma_S M'd}{u\varepsilon}(t_1+2t_2+\varepsilon)\Lambda_h+o(d).
\end{equation}
Here $h(s_d)$ denotes the joint coefficient height of the family $(s_d\otimes\tau)_{\tau\in Z'_d}$ in the fixed monomial basis of $E_d$.
\end{proposition}

\begin{proof}
The generating family $Z'_d$, followed by the two injections in \eqref{eq:injections}, gives the Faltings complex on $\cY$,
\[0\to\cQ_d\xrightarrow{\Phi_d}(\cP^{\otimes d})^{\oplus M'}\xrightarrow{\Psi_d}(\cN_a^{\otimes d(t_1+t_2+\varepsilon)})^{\oplus(M')^2}.\]
Before the last injection, its maps are
\[s\mapsto(s\otimes\tau)_\tau \quad\text{and}\quad (s_\tau)_\tau\mapsto (s_\tau\otimes\tau'-s_{\tau'}\otimes\tau)_{\tau,\tau'}.\]
It is exact at the first two terms: locally a member of $Z'_d$ trivializes $\cR_d$, and the last injection remains injective on the integral variety $\cY$, since it is an isomorphism at $x$. Consequently, with
\[ \cF_d=E_d^{\oplus M'},\qquad n_d=\dim_F\cF_d,\]
the kernel of the induced map on $\cF_d$ is
$\Phi_dH^0(\cY,\cQ_d)$ and has dimension $k_d$.

The injection $\cP^{\otimes d}\hookrightarrow\cN_a^{\otimes dt_1}$, the estimate \eqref{eq:intersection-upper-bound} and the global generation of $\cN_a$ on $\cY$ imply that
\[ n_d\leqslant M'(dt_1)^u\frac{m^u}{u!}\prod_ id_i a_i^{u_i}+o(d^u),\qquad \log n_d=O(\log d).\]
Combining this with \eqref{eq:section-dimension} yields
\begin{equation}\label{eq:dirichlet-ratio}
\frac{n_d}{k_d}\leqslant 4\theta M'(t_1m)^u\prod_i d_i^{1+\omega}+o(1)=\frac{2M'}{u\varepsilon}+o(1),
\end{equation}

We next represent this kernel by linear equations of controlled height. Use the monic relations $Q_j^{(i)}$ to reduce powers of the $W_j^{(i)}$, and the primitive-element relations $P_{u_i+1}^{(i)}$ to express the remaining finitely many monomials in a basis over the rational function field of $\prod_i\bP^{u_i}$. Clearing their denominators once, independently of $d$, gives the construction of \cite[proof of Proposition~4.2]{Remond2005}; if $Y_i=\bP^{N_i}$, no reduction is needed in that factor. Equating the resulting monomial coefficients produces a matrix $A_d$ over the fixed field $F$ with
\begin{equation}\label{eq:system-height}
 \ker A_d=\Phi_dH^0(\cY,\cQ_d),\qquad h(A_d)\leqslant d(t_1+2t_2+\varepsilon)\Lambda_h+o(d).
\end{equation}
The height estimate follows from the local calculation of \cite[Lemme~2.5]{Remond2005}. More explicitly, a coefficient of the input family multiplies a polynomial of the form $\pm\mu P_j^d$, where $\mu$ is a unit monomial of multidegree $d(t_1+\varepsilon)a$. The reduction and the fixed clearing of denominators give
\[\begin{split}
h(A_d)\leqslant{} &h(P_1^d,\dotsc,P_M^d)\\
&+d(t_1+t_2+\varepsilon)\sum_i a_i\bigl(h(B_i)+\beta_S\log(2u_i+2)\bigr)+o(d).
\end{split}\]
The estimate for powers includes the archimedean term noted in \cite[proof of Proposition~3.1]{Dill2020}:
\[ h(P_1^d,\dotsc,P_M^d) \leqslant d\sum_i a_i\delta_i+\beta_S d t_2\sum_i a_i\log(N_i+1).\]
Equation~\eqref{eq:relation-coefficients} and the definition of $\Lambda_h$ imply
\[\sum_i a_i\bigl(h(B_i)+\delta_i+\beta_S\log(2u_i+2)\bigr)\leqslant\Lambda_h,\qquad\beta_S\sum_i a_i\log(N_i+1)\leqslant\Lambda_h.\]
These inequalities prove \eqref{eq:system-height}. The fixed
denominators and polynomial numbers of monomials contribute only
$o(d)$; they do not enlarge the field of definition as $d$ varies.

If $A_d=0$, any coordinate vector gives the assertion. Otherwise choose a maximal independent family of its rows. Its size is $b_d=n_d-k_d$, its kernel is unchanged, and each row has height at most $h(A_d)$. Fix $R>1$ larger than the limiting upper bound in \eqref{eq:dirichlet-ratio}. Hypothesis~\ref{hyp:siegel} therefore supplies a nonzero vector $z_d\in\ker A_d$ satisfying
\[ h(z_d)\leqslant\frac{n_d}{k_d}\bigl(\sigma_S h(A_d)+\gamma_S(F,R)\log n_d+c_S(F,R)\bigr).\]
The corresponding section $s_d$ has $h(s_d)=h(z_d)$ by definition. Now \eqref{eq:dirichlet-ratio} and \eqref{eq:system-height} prove \eqref{eq:small-section-height}; the terms involving $\log n_d$ and $c_S(F,R)$ become $o(d)$. This is the only use of Hypothesis~\ref{hyp:siegel} in the construction.
\end{proof}

For the usual number field height, $\sigma_S=\beta_S=1$ and \eqref{eq:small-section-height} recovers Dill's bound. 

\section{The index is small}\label{sec:small-index}

The product theorem needed here follows from the algebraic argument in \cite{Remond2001Product} and its local height calculations, using the intersection heights of Chen--Moriwaki. The geometric versions in \cite[\S4.5]{Liu2026} and \cite{Jiang2026} have no archimedean error terms. We retain these terms with their factor $\beta_S$. All components below are taken over $\Kbar$. When necessary, we work over a finite extension of $K$ over which all geometric objects are defined, using normalized heights throughout.

Put $\bP=\prod_{i=1}^m\bP^{u_i}$, where $u_i\geqslant 0$ and $u=\sum_i u_i\geqslant 1$. For $y\in\bP(\Kbar)$ and a nonzero multihomogeneous polynomial $G$ of multidegree $e=(e_i)\in(\bN_+)^m$, define the index of $G$ at $y$ with respect to weight $e$ to be 
\[ \ind_{e}(G;y) =\min_{\partial^{[\kappa]}G(y)\neq 0}\sum_i\frac{\abs\kappa_i}{e_i}, \qquad\partial^{[\kappa]}=\frac{1}{\kappa!}\partial^\kappa.\]
Here the derivatives range over all homogeneous coordinates in each block. By the Euler identities this equals the index computed after dehomogenization on any affine chart containing $y$. Write $Z_t(G)$ for the closed locus where this index is greater than $t$, cut out by the divided derivatives of weight at most $t$.

Let
\[ s_n=\frac12\sum_{k=1}^n\sum_{j=1}^k\frac1j\leqslant n\log(n+1)\]
be the Stoll number. Thus $h(\bP^n)=\beta_Ss_n$, with the height convention of \eqref{eq:variety-height}\cite[Proposition~3.7.5]{ChenMoriwaki2021}. Set $\abs e=\sum_ie_i$ and
\[ \cH_e(G)=h(G)+\beta_S\Bigl(\sum_ie_i(s_{u_i}+\log2)+\sum_i\sqrt{u_i}+\frac{u-1}{2}\log\abs e\Bigr).\]

\begin{lemma}\label{lem:common-component}
Let $T$ be a common irreducible component of $Z_t(G)$ and $Z_{t+\eta}(G)$, where $t\geqslant 0$ and $\eta>0$, and put $c=\codim_{\bP}T$. Define
\[ \chi_m(c)=\max_{\substack{\alpha\in\bN^m\\\abs\alpha=c}}\frac{c!}{\alpha!}.\]
If $e_i/e_{i+1}\geqslant 1$ and $e_i/e_{i+1}>\chi_m(c)\eta^{-c}$ for every $i<m$, then $T=\prod_i T_i$ with integral $T_i\subset\bP^{u_i}$. Writing $c_i=\codim T_i$, $\Delta_i=\deg T_i$, and
$R=c!\eta^{-c}/\prod_i c_i!$, one has
\begin{equation}\label{eq:common-component-bounds}
\prod_i\Delta_i\leqslant R,\qquad e_i h(T_i)\prod_{j\neq i}\Delta_j\leqslant c_i R\cH_e(G)\quad(c_i>0).
\end{equation}
\end{lemma}

\begin{proof}
We explain the height input to the algebraic proof of \cite[Proposition~2.1]{Remond2001Product}. For an integral subvariety $T$ of codimension $c$, use the mixed degrees and heights
\[ D_\alpha(T)=\Bigl(\prod_i L_i^{u_i-\alpha_i}\cdot T\Big), \qquad H_\gamma(T)=\Bigl(\prod_i\overline L_i^{u_i-\gamma_i}\mid T\Bigr)_S,\]
where $L_i$ is the pullback of $\cO(1)$ from the $i$-th factor, with its standard metric, $\alpha,\gamma\in\bN^m$ with $\abs\alpha=c$, $\abs\gamma=c-1$ and $\alpha_i,\gamma_i\leqslant u_i$. These heights are nonnegative, by coordinate-section induction as in Lemma~\ref{lem:bezout}, and are additive on effective cycles.

Let $I$ be generated by a finite family $\cA$ of polynomials of multidegrees at most $e$, and let $T$ be a component of its zero locus, with codimension $c$. The construction of \cite[Proposition~3.2 and \S5]{Remond2001Product} gives $c$ elements of $I$ whose intersection has finite length $\lambda$ at the generic point of $T$, with
\begin{align}
\lambda D_\alpha(T)&\leqslant\frac{c!}{\alpha!}e^\alpha, \label{eq:product-degree-intersection}\\
\lambda H_\gamma(T)&\leqslant\frac{c!}{\gamma!}e^\gamma C+\beta_S\sum_i s_{u_i}\frac{c!}{(\gamma+\boldsymbol{\varepsilon}_i)!}e^{\gamma+\boldsymbol{\varepsilon}_i}, \label{eq:product-height-intersection}
\end{align}
where
\[ C=h(\cA)+\beta_S\Bigl(\frac{c-1}{2}\log\abs e+\sum_i\sqrt{u_i}\Bigr).\]
The mixed B\'ezout inequalities used in these estimates follow as in the proof of Lemma~\ref{lem:bezout}. Apply \cite[Proposition~4.4.4(3)]{ChenMoriwaki2021} to the standard bundles pulled back from the factors: the resulting positive mixed intersection measure has mass equal to the corresponding mixed degree. Bounding the section norm in each block as in that proof and using multilinearity give the required inequalities. The bounded integer linear combinations in \cite[\S5]{Remond2001Product} are unchanged. The combinatorial errors occur only at archimedean places, so their integrated bounds acquire the factor $\beta_S$; the initial mixed heights of $\bP$ give the Stoll terms in \eqref{eq:product-height-intersection}. 

Thus these are adelic versions of the cited complete-intersection estimates. In the geometric case they reduce to the estimates of \cite[Lemma~4.9]{Liu2026} and the product theorem appendix in \cite{Jiang2026}.

Apply this to the derivatives of $G$ of weight at most $t$. Their joint coefficient height is at most $h(G)+\beta_S\abs e\log2$. The common-component condition implies that derivatives of elements in $I$ of weight at most $\eta$ vanish on $T$. The algebraic multiplicity estimate and length comparison of \cite[Proposition~3.1 et Lemme~3.1]{Remond2001Product} therefore give
\[ \lambda\geqslant\eta^c\prod_ie_i^{v_i},\]
where
\[ v_i=\codim T^{(i)}-\codim T^{(i+1)},\]
with $T^{(i)}\subset\prod_{j=i}^m\bP^{u_j}$ denoting the image of $T$ under projection. Each codimension is taken in the corresponding product, with the convention $\codim T^{(m+1)}=0$. Comparing this with \eqref{eq:product-degree-intersection}, the argument in \cite[proof of Proposition~2.1]{Remond2001Product} forces $T$ to be a product under the assumption on $e$. This algebraic argument works in characteristic zero and applies over $\Kbar$.

For this product $v_i=c_i$ and, with $\mathbf c=(c_i)$,
\[ D_{\mathbf c}(T)=\prod_i\Delta_i, \qquad H_{\mathbf c-\boldsymbol{\varepsilon}_i}(T)=h(T_i)\prod_{j\neq i}\Delta_j\quad(c_i>0).\]
To see the height identity, put $r_j=\dim T_j=u_j-c_j$. For fixed $i$ and $j\neq i$, cutting $T_j$ with $r_j$ general hyperplanes gives a zero cycle of degree $\Delta_j$. At each cut, the local integral term in \cite[Proposition~4.4.4(3)]{ChenMoriwaki2021} vanishes: its mixed measure contains $r_i+1$ Chern factors pulled back from the $r_i$-dimensional variety $T_i$, and hence has total mass $0$. The remaining cycle projects finitely onto $T_i$ with degree $\prod_{j\neq i}\Delta_j$, counting
multiplicities, and the remaining metrized line bundles are pulled back from $\overline{\mathcal O(1)}$ on $T_i$. The arithmetic projection formula \cite[Theorem~4.4.9]{ChenMoriwaki2021} therefore gives the height identity. Use $\alpha=\mathbf c$ and $\gamma=\mathbf c-\boldsymbol{\varepsilon}_i$ in the preceding estimates, divide by $\eta^ce^{\mathbf c}$, and use
$c\leqslant u$ and $(\gamma+\boldsymbol{\varepsilon}_j)!\geqslant\gamma!$.
This gives \eqref{eq:common-component-bounds}.
\end{proof}

\begin{theorem}\label{thm:adelic-product}
Let $0<\tau/m\leqslant(2u^2)^{-u}$ and put $L=(m/\tau)^u$. Suppose that $e_i/e_{i+1}\geqslant L$ for every $i<m$ and that $\ind_e(G;y)>\tau$.
Then there exist $l$ with $u_l\geqslant 1$ and a nonzero homogeneous
polynomial $U$ in the $l$-th block such that
\begin{align}
U(y_l)&=0,\qquad 1\leqslant\deg U<L,\label{eq:product-equation-degree}\\
e_l h(U)&\leqslant u_l L\cH_e(G)+2\beta_Se_l u_l L\log\bigl(2L(u_l+1)\bigr).\label{eq:product-equation-height}
\end{align}
\end{theorem}

\begin{proof}
Write 
\[ t_j=\frac{j\tau}{2u^2}\quad(0\leqslant j<u),\qquad t_u=\tau.\]
Note that $\chi_m(j)\leqslant m^{j-1}$ for $j\geqslant 1$, therefore
\begin{equation}\label{eq:product-index-levels}
\chi_m(j)(t_j-t_{j-1})^{-j}<L\qquad(1\leqslant j\leqslant u).
\end{equation}
Indeed, for $j<u$ the ratio of the left side to $L$ is at most $(2u^2)^j(\tau/m)^{u-j}\leqslant (2u^2)^{-1}$; for $j=u$ it is at most
\[m^{-1}\Bigl(1-\frac{u-1}{2u^2}\Big)^{-u}\leqslant m^{-1}\Bigl(1-\frac{u-1}{2u}\Bigr)^{-1}<1.\]

Choose nested components $T^{(u)}\subset\dotsb\subset T^{(0)}$ through $y$, where $T^{(j)}$ is a component of $Z_{t_j}(G)$. Their codimensions lie between $1$ and $u$. For the first $j\geqslant 1$ with $\codim T^{(j)}\leqslant j$, one has $T^{(j)}=T^{(j-1)}$ and both have codimension $j$. Lemma~\ref{lem:common-component} and \eqref{eq:product-index-levels} then give a proper product $T^{(j)}=\prod_i T_i$ through $y$ with
\[ \prod_i\deg T_i<L,\qquad e_i h(T_i)\leqslant u_iL\cH_e(G)\quad(T_i\neq\bP^{u_i}).\]
Choose a proper factor $T_l$ and let $e=\deg T_l<L$. Apply the adapted-coordinate construction of section~\ref{sec:reduction} to $T_l\subset\bP^{u_l}$. The coefficient bound \eqref{eq:relation-coefficients}, with ambient dimension $u_l$, variety dimension $r$, and degree $e$, gives a nonzero degree-$e$ relation $F$ in the first $r+2$ adapted coordinates with
\[ h(F)\leqslant h(T_l)+\beta_S e(r+1)\log\bigl(e(u_l+1)\bigr).\]
The estimate \eqref{eq:coordinate-height-change} bounds the height increase on returning to the original coordinates. Using $r+1\leqslant u_l$, we obtain a nonzero degree-$e$ equation $U$ of
$T_l$ satisfying
\[\begin{split}
h(U)&\leqslant h(T_l)+\beta_S\bigl(e u_l\log\bigl(e(u_l+1)\bigl)+e\log\bigl(e u_l(u_l+1)\bigr)\bigr)\\
&\leqslant h(T_l)+\beta_S\bigl(e(u_l+1)\log\bigl(2e(u_l+1)\bigr)\bigr)\\
&\leqslant h(T_l)+2\beta_S u_l L\log\bigl(2L(u_l+1)\bigr). 
\end{split}\]
This proves \eqref{eq:product-equation-height}.
\end{proof}

When $\beta_S=0$ and $e_i=db_i$, the conclusion reads $b_l h(U)\leqslant u_l L h(G)/d$, as in the product theorem appendix of \cite{Jiang2026}; see also \cite[Lemmas~4.10--4.13]{Liu2026}. 

Return to the minimal product and the section $s_d$ of Proposition~\ref{prop:small-section}. On a neighborhood of $x$ in $Y\cap U^0$, trivialize $\cQ_d$ and use the local coordinates $V_j^{(i)}/V_0^{(i)}$ furnished by Corollary~\ref{cor:local-coordinates}. Let $\iota_d$ be the index of $s_d$ at $x$ for the weights $dt_1a$. It is independent of the trivialization. Put $D=\prod_i d_i$ and
$\rho=\prod_i\rho_i:Y\to\bP$.

\begin{lemma}[Norm polynomial]
\label{lem:norm-polynomial}
There is a nonzero multihomogeneous polynomial $G_d$ of multidegree
$dt_1Da$ such that
\begin{align}
\ind_{dt_1Da}(G_d;\rho(x))&\geqslant\iota_d/D,\label{eq:norm-index}\\
h(G_d)&\leqslant D\bigl(h(s_d)+dt_1\sum_i a_i\bigl(h(B_i)+\beta_S\log(2u_i+2)\bigr)\bigr)+o(d). \label{eq:norm-height}
\end{align}
\end{lemma}

\begin{proof}
Choose $\tau\in Z'_d$ nonzero at $x$ and write the image of $s_d\otimes\tau$ under the first injection of \eqref{eq:injections} as $\alpha\prod_i(V_0^{(i)})^{dt_1a_i}$. Let $N(\alpha)$ be the norm of $\alpha$ for $\Kbar(Y)/\Kbar(\bP)$, an extension of degree $D$, and set
\[ G_d=\prod_i(V_0^{(i)})^{dt_1Da_i}N(\alpha).\]
The integrality argument in \cite[proof of Lemma~4.1]{Dill2020} shows that this expression is a polynomial: the corresponding norm on the multihomogeneous coordinate ring is integral over the polynomial ring in the $V_j^{(i)}$, $j\leqslant u_i$, which is integrally closed. Its multidegree is $dt_1Da$. Since $\alpha\in\cO_{Y,x}$ is integral over $\cO_{\bP,\rho(x)}$, $N(\alpha)/\alpha$ is also integral over $\cO_{\bP,\rho(x)}$, thus integral over $\cO_{Y,x}$. As $\rho$ is \'etale at $x$, $\cO_{Y,x}$ is normal and then $N(\alpha)/\alpha\in\cO_{Y,x}$ is regular at $x$. This proves \eqref{eq:norm-index}.

For the height, use the coefficient calculation of \cite[Lemme~5.5]{Remond2005}, based on its Lemme~2.5. The archimedean error in the Lemme~2.5 contributes
\[\beta_Sdt_1\sum_i a_i\log(2u_i+2)+o(d)\]
to the height. Taking the norm multiplies the resulting height bound by $D$. The fixed coefficients in the norm expansion contribute $O(1)$, while the numbers of monomials grow polynomially in $d$ and hence contribute only $O(\log d)$. These latter contributions are absorbed in $o(d)$, giving \eqref{eq:norm-height}. For the geometric version, see also \cite[Lemma~4.12]{Liu2026}.
\end{proof}

\begin{proposition}\label{prop:small-index}
For sufficiently large $d\in d_0\bN$, the index $\iota_d$ of the section constructed in Proposition~\ref{prop:small-section} satisfies
\begin{equation}\label{eq:small-index}
\iota_d\leqslant\frac{\varepsilon}{4t_1\max_i d_i(N_i+1)}.
\end{equation}
\end{proposition}

\begin{proof}
Suppose that \eqref{eq:small-index} fails, and put
\[ q=\psi(u),\qquad \tau=m\Lambda^{-\omega'q},\qquad L=(m/\tau)^u=\Lambda^{\omega'u q}.\]
Recall $A=\max_i(N_i+1)$ and $\sD=\prod_i\deg X_i$. By \eqref{eq:minimal-product} and  \eqref{eq:epsilon},
\[ \frac{\iota_d}{D}>\frac{D^{-\omega'}}{8u\theta t_1^{u+1}m^u A}\geqslant\frac{\sD^{-\omega'}\Lambda^{-\omega'(q-1)}}{8u\theta t_1^{u+1}m^u A}\geqslant\tau.\]
The last inequality uses $\Lambda^{\omega'}\geqslant8u\theta t_1^{u+1}m^{u+1}A\sD^{\omega'}$, which follows from the definition of $\Lambda$ \eqref{eq:Lambda}. Also $a_i/a_{i+1}\geqslant\Lambda^{\psi(0)}\geqslant L$, since $\psi(0)\geqslant\omega'u\psi(u)$, and $\tau/m\leqslant(2u^2)^{-u}$. Thus Theorem~\ref{thm:adelic-product} applies to $G_d$ at $\rho(x)$. It gives $l$ and $U$ with
\[ U(\rho_l(x_l))=0,\qquad 1\leqslant\deg U<L.\]

Divide \eqref{eq:product-equation-height} by $dt_1D$. Since the terms $\sum_i\sqrt{u_i}$ and $\log\abs{dt_1Da}$ are $o(d)$, \eqref{eq:norm-height} gives
\[\begin{split}
a_lh(U)\leqslant u_lL\biggl(&\frac{h(s_d)}{dt_1}+\sum_i a_i\bigl(h(B_i)+\beta_S\log(2u_i+2)\bigr)\\
&+\beta_S\sum_i a_i(s_{u_i}+\log2)+2\beta_Sa_l\log\bigl(2L(u_l+1)\bigr)\biggr)+o(1).
\end{split}\]
By \eqref{eq:Lambda-h} and \eqref{eq:relation-coefficients}, each of the two sums is at most $\Lambda_h$; also $\beta_Sa_l\log(2u_l+2)\leqslant\Lambda_h$. Using \eqref{eq:small-section-height}, $t_1\geqslant 1$, and $\sigma_SM'\geqslant 1+\sum_i N_i> u$, we obtain
\begin{equation}\label{eq:factor-equation-budget}
a_lh(U)\leqslant u_lL\Bigl(\frac{2\sigma_SM'}{u\varepsilon}(1+2t_2+3\varepsilon)\Lambda_h+2\beta_Sa_l\log L\Bigr)+o(1).
\end{equation}

Equations
\eqref{eq:Lambda} and~\eqref{eq:minimal-product} give
\[ 4u\theta A^m(t_1m)^u\leqslant\Lambda^m, \qquad D^{1+\omega}\leqslant\Lambda^{(1+\omega)q}.\]
Since $m\leqslant u\psi(u)$ and $\omega'\geqslant 2$, it follows that
\[\frac{2\sigma_S M' u_l}{u\varepsilon}\leqslant\frac{2\sigma_SM'}{\varepsilon}\leqslant \sigma_SM\Lambda^{m+(1+\omega)q}\leqslant \sigma_SM\Lambda^{(\omega'u-1)q}.\]
Furthermore $\beta_S a_l u_l\leqslant\Lambda_h$ and $\Lambda\geqslant 32$, so
\[ 2\beta_S a_l u_l\log L \leqslant 2\Lambda_h\omega'u q\log\Lambda\leqslant\Lambda_h\Lambda^{(\omega'u-1)q}.\]
Equation~\eqref{eq:factor-equation-budget} therefore implies
\[ a_l h(U)\leqslant \sigma_SM\Lambda^{(2\omega'u-1)q} (2+2t_2+3\varepsilon)\Lambda_h+o(1).\]
Now $\varepsilon\leqslant 1/4$, hence $2+2t_2+3\varepsilon\leqslant 5t_2$. Also $d_l\leqslant\deg(X_l)\Lambda^{q-1}$ and $\Lambda\geqslant 32\deg(X_l)$, so $20d_l<\Lambda^q$. For sufficiently large $d$, the last bound is consequently strictly less than
\[ \Lambda^{2\omega'u q}\frac{B}{4d_l}\Lambda_h=\Lambda^{2\psi(u-1)-2\psi(u)}\frac{B}{4d_l}\Lambda_h.\]
We have obtained all three conditions forbidden by Proposition~\ref{prop:no-small-equation}, a contradiction.
\end{proof}

\section{Finishing the proof}\label{sec:finishing}

We finish with the derivative argument of \cite[\S5]{Dill2020}. There are two archimedean corrections to keep track of: the term for powers of $P_i$, already included in \eqref{eq:system-height}, and the derivative constant in \cite[proof of Proposition~5.3]{Remond2005}, where Dill replaces $8$ by $12$. Both acquire the factor $\beta_S$ upon integration. They disappear for the geometric heights in \cite{Jiang2026}.

On a neighborhood of $x_i$, put
\[ v_j^{(i)}=V_j^{(i)}/V_0^{(i)}\quad(1\leqslant j\leqslant u_i), \qquad w_j^{(i)}=W_j^{(i)}/V_0^{(i)}\quad(0\leqslant j\leqslant N_i).\]
By Corollary~\ref{cor:local-coordinates}, the $v_j^{(i)}$ are \'etale coordinates. For $\ell\in\bN^{u_i}$, let
\[ \partial_i^{[\ell]}=\prod_{j=1}^{u_i}\frac{1}{\ell_j!}\Bigl(\frac{\partial}{\partial v_j^{(i)}}\Bigr)^{\ell_j}.\]
Write $J_i=\{j:W_j^{(i)}\lvert_{Y_i}\not\equiv0\}$ and choose
\[ \gamma_i=\prod_{j\in J_i}\Bigl(\frac{w_j^{(i)}(x_i)^{b_{i,j}}}{b_{i,j}!}\frac{\partial^{b_{i,j}}Q_j^{(i)}}{\partial T_{u_i+1}^{b_{i,j}}}(1,v_1^{(i)}(x_i),\ldots,v_{u_i}^{(i)}(x_i),w_j^{(i)}(x_i))\Bigr)^{1/b_{i,j}}.\]
These are nonzero by the same corollary. Enlarge the fixed field of definition $F$ finite over $K$ to contain $x$ and the chosen roots $\gamma_i$; this enlargement is independent of $d$. For $v\in\Omega_F$, put $\varepsilon_v=1$ if $v$ is archimedean and $\varepsilon_v=0$ otherwise.

\begin{lemma}[Local derivative estimate]\label{lem:local-derivatives}
For every $v\in\Omega_F$, $\ell\in\bN^{u_i}$, and
$0\leqslant j\leqslant N_i$, one has
\[ \abs{\partial_i^{[\ell]}w_j^{(i)}(x_i)}_v \leqslant\abs{w_j^{(i)}(x_i)}_v\bigl(\abs{\gamma_i}_v^{-2}C_{i,v}\bigr)^{\abs\ell},\]
where
\[ C_{i,v}=2^{-\varepsilon_v}\Bigl(\bigl(d_i(N_i+1)\bigr)^{6d_i\varepsilon_v}\max_{b\in B_i}\abs b_v\max_{0\leqslant k\leqslant N_i}\abs{w_k^{(i)}(x_i)}_v^{d_i}\Bigr)^{2(N_i+1)}.\]
\end{lemma}

\begin{proof}
The proof of \cite[Lemma~5.1]{Dill2020} works verbatim: it uses only the monic relations $Q_j^{(i)}$, the nonvanishing just established, and the local absolute value.
\end{proof}

Choose $d\in d_0\bN$ sufficiently large for Propositions~\ref{prop:small-section} and~\ref{prop:small-index}, and choose $\kappa=(\kappa_i)_i$ realizing the index $\iota_d$ of $s_d$. Thus,
\begin{equation}\label{eq:derivative-orders}
\sum_i\frac{\abs{\kappa_i}}{dt_1a_i}=\iota_d,\qquad\abs{\kappa_i}\leqslant\frac{d\varepsilon a_i}{4d_i(N_i+1)}.
\end{equation}
In a local trivialization of $\mathcal Q_d$, the operator $\sD_\kappa=\prod_i\partial_i^{[\kappa_i]}$ is nonzero on $s_d$ at $x$, whereas every derivative of smaller weight vanishes there.

The first injection in \eqref{eq:injections} and the nonvanishing of $V_0^{(i)}$ at $x$ give a local frame for $\cP$. In this frame a coordinate $\Xi_j$ is represented by a monomial $\xi_j$ of multidegree $t_1a$ in $w_k^{(i)}$. If $\xi^\nu$ ranges over the degree-$d$ monomials corresponding to the basis of $E_d$, the divided Leibniz formula and Lemma~\ref{lem:local-derivatives} give
\begin{equation}\label{eq:monomial-derivatives}
\abs{\sD_\kappa(\xi^\nu)(x)}_v\leqslant \abs{\xi^\nu(x)}_v\prod_i2^{\varepsilon_v(\abs{\kappa_i}+d t_1u_i a_i)}\bigl(\abs{\gamma_i}_v^{-2}C_{i,v}\bigr)^{\abs{\kappa_i}}.
\end{equation}
Indeed, the number of terms in the $i$-th block is
\[\prod_{j=1}^{u_i}\binom{\kappa_{i,j}+dt_1a_i-1}{dt_1a_i-1}
\leqslant 2^{\abs{\kappa_i}+dt_1u_i a_i},\]
see \cite[Proof of Proposition~5.3]{Remond2005}. 

We next explain the passage from these local estimates to the height of the section. Keeping the above frame for $\cP^{\otimes d}$, choose compatible local frames for $\mathcal Q_d$ and $\mathcal R_d$, and write $s_d=f e_Q$ and $\tau=g_\tau e_R$ for $\tau\in Z'_d$. Minimality of the derivative weight gives
\[ \sD_\kappa(fg_\tau)(x)=\sD_\kappa(f)(x)g_\tau(x),\qquad \sD_\kappa(f)(x)\neq 0.\]
Thus the vector of these derivatives represents the same projective point as $Z'_d(x)$. By the product formula and the coefficient definition of $h(s_d)$,
\begin{equation}\label{eq:differentiated-section-height}
\begin{split}
h(Z'_d(x)) &=h\bigl((\sD_\kappa(fg_\tau)(x))_\tau\bigr)\\
&\leqslant h(s_d)+h\bigl((\sD_\kappa(\xi^\nu)(x))_\nu\bigr)+\beta_S\log\dim E_d.
\end{split}
\end{equation}
The last term is $o(d)$. 

Integrate \eqref{eq:monomial-derivatives}, using $\max_\nu\abs{\xi^\nu(x)}_v\leqslant (\max_j\abs{\xi_j(x)}_v)^d$ and the product formula for $\gamma_i$. The factor $2^{-\varepsilon_v}$ in $C_{i,v}$ cancels the contribution $\abs{\kappa_i}\varepsilon_v\log2$ and we obtain
\[\begin{split}
h\bigl((\sD_\kappa(\xi^\nu)(x))_\nu\bigr)\leqslant {}&d h(\Xi(x))+\beta_Sdt_1\sum_i u_i a_i\log2\\
&+\sum_i2(N_i+1)\abs{\kappa_i}\bigl(h(B_i)+d_i h(x_i)\bigr)\\
&+12\beta_S\sum_i d_i(N_i+1)\abs{\kappa_i}\log\bigl(d_i(N_i+1)\bigr).
\end{split}\]
Write $H=h_{\cN_a}(x)$. By \eqref{eq:derivative-orders}, the point-height terms in this bound are at most $d\varepsilon H/2$. Moreover, \eqref{eq:relation-coefficients} and \eqref{eq:Lambda-h} imply
\[\begin{split}
\sum_i&2(N_i+1)\abs{\kappa_i}\bigl(h(B_i)+6\beta_Sd_i\log(d_i(N_i+1))\bigr)\\
&\leqslant d\varepsilon\sum_i a_i\Bigl(\frac{h_i+\beta_Sd_i(u_i+1)\log(d_i(N_i+1))}{2d_i}+3\beta_S\log(d_i(N_i+1))\Bigr)\\
&\leqslant 4d\varepsilon\Lambda_h,
\end{split}\]
and $\beta_S\sum_i u_i a_i\log2\leqslant\Lambda_h$.
Finally, the definition of $Z'_d$ gives the exact identity
\[ h(Z'_d(x))-d h(\Xi(x))=d\varepsilon H-d h_{\cM}(x).\]
Combining these facts with \eqref{eq:differentiated-section-height} yields
\begin{equation}\label{eq:final-height-comparison}
\frac{\varepsilon}{2}H-h_{\cM}(x)\leqslant\frac{h(s_d)}{d}+(t_1+4\varepsilon)\Lambda_h+o(1).
\end{equation}

\begin{proof}[Completion of the proof of Theorem~\ref{thm:vojta-remond}]
Since $2\sigma_SM'/(u\varepsilon)\geqslant 1$, Proposition~\ref{prop:small-section} and \eqref{eq:final-height-comparison}, followed by $d\to\infty$, give
\begin{equation}\label{eq:final-section-budget}
\frac{\varepsilon}{2}H-h_{\cM}(x)\leqslant\frac{2\sigma_SM'}{u\varepsilon}(2t_1+2t_2+5\varepsilon)\Lambda_h.
\end{equation}
We bound this coefficient in terms of $\Lambda$ and $B$. Put $q=\psi(u)$ and $p=\max\{3,m\}$. As $u\leqslant u_0$, $m\leqslant u_0$, and $u_0\geqslant 2$, the definitions of $\Lambda$ and $p$ imply
\[ \sigma_St_1 A^m(2u\theta)^2(t_1m)^{2u}\leqslant\sigma_S\Lambda^p.\]
Use $D=\prod_i d_i\leqslant\sD\Lambda^{q-1}\leqslant\Lambda^q$ to deduce
\[\begin{split}
\sigma_S M' t_1\varepsilon^{-2}&=\sigma_SMt_1\prod_i(N_i+1)(2u\theta)^2(t_1m)^{2u}D^{2(1+\omega)}\\
&\leqslant \sigma_SM\Lambda^{p+2(1+\omega)q}\leqslant \sigma_SM\Lambda^{2\omega'u q-2}.
\end{split}\]
For the last step, use $m\leqslant uq$, $uq\geqslant 2$, and $\omega'\geqslant 2$: these give $p\leqslant 4uq-2$ and $p+2(1+\omega)q\leqslant 2\omega'u q-2$. Also $\varepsilon\leqslant 1/4$, so $8(2+2t_2+5\varepsilon)\leqslant 42t_2\leqslant\Lambda^2t_2$. Consequently
\[ \frac{8\sigma_SM'}{u\varepsilon^2}(2t_1+2t_2+5\varepsilon)\leqslant \sigma_SMt_2\Lambda^{2\omega'u q}=B\Lambda^{2\psi(u-1)-2q}.\]
Apply this to \eqref{eq:final-section-budget}, followed by \eqref{eq:Lambda-h-bound}, to obtain
\[ \frac{\varepsilon}{2}H-h_{\cM}(x)\leqslant\frac{\varepsilon}{4}B^{1-u}\Lambda^{2\psi(u-1)-2\psi(0)}\sum_i a_i c_3^{(i)}\leqslant\frac{\varepsilon}{4}H.\]
Here $u\geqslant 1$, $\psi(u-1)\leqslant\psi(0)$, and the last inequality uses \eqref{eq:height-hypotheses}. Thus $H\leqslant 4\varepsilon^{-1}h_{\cM}(x)$. In particular $h_{\cM}(x)>0$, since $H\geqslant\sum_i a_i c_3^{(i)}>0$. Finally,
\[ 4\varepsilon^{-1}=8u\theta(t_1m)^uD^{1+\omega}\leqslant\Lambda^{(1+\omega)q+2}\leqslant\Lambda^{\omega'u q}\leqslant\Lambda^{\psi(0)}=c_1.\]
For the first inequality use $8u\theta(t_1m)^u\leqslant\Lambda^2$ and the same bound for $D$; the other two follow from $\omega'=3+\omega$ and the definition of $\psi$. This proves \eqref{eq:main-inequality}.
\end{proof}

\appendix

\section{A Siegel estimate over arithmetic function fields}\label{app:arithmetic-siegel}

Let $K$ be a finitely generated field over $\bQ$ with $r=\trdeg_\bQ K>0$. This appendix verifies Hypothesis~\ref{hyp:siegel} for the adelic curve structure over $K$ induced by an arithmetic model with nef and big polarizations, a case of \cite[\S3.2.6]{ChenMoriwaki2020}. The proof first treats a diagonal ample polarization; the general case follows by height comparison.

Let $\cB$ be a normal integral scheme, projective and flat over $\Spec\bZ$, with function field $K$, and put $D=r+1$. We first consider the diagonal polarization $\overline H^{\,r}$, where $\overline H$ is the pullback of $\overline{\cO(1)}(1)$ via a projective embedding $\cB\hookrightarrow\bP^N_\bZ$; here $\overline L(t)$ means that the metric is multiplied by $e^{-t}$. Then $\overline H$ is an arithmetically ample smooth hermitian line bundle on $\cB$, in the convention of \cite[p.~608]{Yuan2008}, and $H_\bQ$ is globally generated. In addition, $\overline H(-1)$ is nef. These data define an adelic curve structure $S$ on $K$ as in \cite[\S3.2.6]{ChenMoriwaki2020}.  Write
\[ V=\widehat{\deg}(\overline H^{D})>0, \qquad v=\deg(H_{\bQ}^{r})>0.\]

We use throughout the arithmetic projection formula: for a dominant generically finite morphism $f:\cB'\to\cB$ and integrable metrized line bundles $\overline L_1,\dotsc,\overline L_D$ on $\cB$,
\[\widehat{\deg}_{\cB'}(f^*\overline L_1\dotsm f^*\overline L_D)=\deg(f)\,\widehat{\deg}_{\cB}(\overline L_1\dotsm\overline L_D);\]
see \cite[Theorem~4.4.9]{ChenMoriwaki2021}. In particular, these intersection numbers are unchanged under birational pullback.

\begin{lemma}\label{lem:arithmetic-target}
There is a constant $C$, depending only on $(\mathcal B,\overline H)$, with the following property. Let $\pi:\mathcal B'\to\mathcal B$ be a normal projective birational model, and let $\overline L$ be a nef line bundle with smooth hermitian metric on $\mathcal B'$. Set
\[ q=\frac{\widehat{\deg}((\pi^*\overline H)^r\overline L)}{V}.\]
For every integer $p\geqslant\max\{2,q\}$, one has
\begin{equation}\label{eq:uniform-arithmetic-target}
\widehat h^0(p\pi^*\overline H+\overline L)\leqslant\frac{V}{D!}(p+q)^D+C p^{D-1}\log p,
\end{equation}
where $\widehat h^0$ is the logarithm of the number of integral sections of supremum norm at most $1$. In particular, $C$ is independent of $\pi$ and $\overline L$.
\end{lemma}

\begin{proof}
Suppress pullbacks in the notation and put $\overline T=p\overline H+\overline L$. Then we have 
\begin{equation}\label{eq:arithmetic-target-volume}
\widehat{\vol}(\overline T)=\widehat{\deg}(\overline T^{\,D})\leqslant\frac{\widehat{\deg}(\overline H^{\,r}\overline T)^D}{V^{D-1}}=V(p+q)^D.
\end{equation}
The first equality follows from \cite[Corollary~5.5]{Moriwaki2009}. The inequality is the arithmetic Khovanskii--Teissier inequality for nef hermitian line bundles \cite[Theorem~2.9(2)]{Ikoma2015}. It remains to control the difference between $\widehat h^0$ and volume uniformly.

Let $k$ be the algebraic closure of $\bQ$ in $K$. Normality implies that $\cB'$ is an $\cO_k$-scheme, where $\cO_k$ is the ring of integers in $k$. The generic fibre $\cB'\otimes_{\cO_k}k=\cB'\otimes_\bZ\bQ$ is geometrically integral as a $k$-scheme. Put $\delta=\deg(H_{\bQ}^{r-1}L_{\bQ})$. Nonnegativity of arithmetic intersections of nef bundles \cite[Proposition~2.3(1)]{Moriwaki2000} implies
\begin{equation}\label{eq:generic-degree-from-height}
Vq=\widehat{\deg}(\overline H^r\overline L)=\widehat{\deg}\bigl(\bigl(\overline H(-1)\bigr)^r\overline L\bigr)+r\deg(H_\bQ^{r-1}L_\bQ)\geqslant r\delta,
\end{equation}
by expanding $\overline H=\overline H(-1)+\overline{\mathcal O}(1)$. Thus $\delta=O(p)$, uniformly on the birational model $\cB'$. For $M=\pi^*H_k$, the Siu inequality \cite[Theorem~2.2.15]{Lazarsfeld2004} shows that $sM-L_k=(p+s)M-T_k$ is big for
\[s=\floor{r\delta_k/v_k}+1=\floor{r\delta/v}+1=O(p),\]
where
\[v_k=\deg(H_{k}^{r})=v/[k:\bQ],\qquad \delta_k=\deg(H_k^{r-1}L_k)=\delta/[k:\bQ].\]
Thus 
\[\deg(T_k^r)\leqslant (p+s)^r\deg(H_k^r)=O(p^r). \]
On the other hand, as $M$ is globally generated, the invariant $\varepsilon(T_k)$ used in \cite{YuanZhang2018} before its Theorem~1.1 satisfies
\[ \varepsilon(T_k)\leqslant (p+s+1)^{r-1}\deg(H_k^r)=O(p^{r-1}).\]
Together with \cite[Theorem~1.1]{YuanZhang2018}, this shows that
\[\rho:=h^0(\cB'_{\bQ},T_{\bQ})=[k:\bQ]h^0(\cB'_k,T_k)=O(p^r).\]
All constants in these estimates depend only on the fixed polarized model. With notation in \cite[p.~259]{YuanZhang2018}, we have 
\[ \dvol(\overline T)=\deg(T_k^r)\geqslant p^r\deg(H_k^r)\]
because $\overline T$ is nef, and both $H_k$ and $L_k$ are nef.

Since 
\[\widehat{\vol}(\overline T)=\widehat{\deg}((p\overline H+\overline L)^D)\geqslant p^D V>0,\]
$\overline T$ is big and the effective arithmetic estimate \cite[Theorem~1.3]{YuanZhang2018} gives
\[ \widehat h^0(\overline T)\leqslant\Bigl(\frac1{D!}+\frac{r\varepsilon(T_k)}{\dvol(\overline T)}\Bigr)\widehat{\vol}(\overline T)+4\rho\log(3\rho).\]
Its first error term is $O(p^{D-1})$ by \eqref{eq:arithmetic-target-volume}, and its second is $O(p^{D-1}\log p)$. This proves \eqref{eq:uniform-arithmetic-target} with the claimed uniformity. The proof uses only the assumptions that $\overline H(-1)$ is nef, $V>0$, and $H_k$ is big and globally generated. Thus the same bound holds for a finite pullback of $\overline H$, with a constant depending on that fixed finite model.
\end{proof}

Let $F$ be a finite extension of $K$ and $w=[w_1:\dotsb:w_n]\in\bP^{n-1}(F)$. The point w defines a rational map $\cB_F\to\bP^{n-1}_\bZ$, which extends to a morphism $f_{w}$ after a birational modification $\cB'\to\cB_F$. Let $(X_j)$ be coordinates on $\bP^{n-1}$, and 
\[ \overline E_{w,\max} =f_w^*(\cO(1),\norm\cdot_{\max}), \]
where
\[\norm{X_j}_{\max}(x) =\frac{\abs{x_j}}{\max_\ell\abs{x_\ell}},\qquad x=[x_1:\dotsb:x_n]\in\bP^{n-1}(\bC). \]
Then $\overline E_{w,\max}$ is nef.

By \cite[Proposition~3.3.2]{Moriwaki2000},
\[ h_S(w)=\widehat{\deg}(\overline H_F^r\overline E_{w,\max})/[F:K].\]
Evaluating the height of $[1:2]$ both adelically and intersection-theoretically gives
\[\beta_S\log2=\frac{1}{[F:K]}\widehat{\deg}(\overline H_F^{r}\overline{\cO}(\log2))=\frac{1}{[F:K]}\deg(H_{F,\bQ}^r)\log2\]
and
\[\beta_S=\deg(H_{F,\bQ}^r)/[F:K].\]
Also recall that the standard metric on $\overline{\cO(1)}$ used throughout is given by the hermitian metric
\[\norm{X_j}_{2}(x) =\frac{\abs{x_j}}{\sqrt{\abs{x_1}^2+\dotsb+\abs{x_n}^2}}.\]
Since $\norm{X_j}_{\max}/\sqrt n\leqslant\norm{X_j}_2\leqslant\norm{X_j}_{\max}$, we have 
\[h_S(w)\leqslant\widehat{\deg}(\overline H_F^rf_w^*\overline{\cO(1)})/[F:K]\leqslant h_S(w)+\frac12\beta_S\log n.\]

\begin{proposition}[Arithmetic Siegel estimate]
\label{prop:arithmetic-siegel}
Fix the adelic curve structure on $K$ induced by the diagonal ample model polarization above. For every finite extension $F/K$ and every $R>1$, there is a constant $c(F,R)$ such that every matrix $A\in M_{b\times n}(F)$ of rank $b$, with $1\leqslant b<n$ and $\lambda:=n/(n-b)\leqslant R$, has a nonzero solution $z\in F^n$ of $Az=0$ satisfying
\begin{equation}\label{eq:arithmetic-siegel}
  h(z)\leqslant D\frac{n}{n-b}\bigl(\max_i h(A_i)+\beta_S\log n\bigr)+c(F,R).
\end{equation}
Thus Hypothesis~\ref{hyp:siegel} holds with
$\sigma_S=\trdeg_{\bQ}K+1$.
\end{proposition}

\begin{proof}
Let $\cB_F$ be the normalization of $\cB$ in $F$, and let $\overline H_F$ be the pullback of $\overline H$. It is nef and big, $\overline H_F(-1)$ is nef, and $V_F:=\widehat{\deg}(\overline H_F^D)=[F:K]V$.

There exists a normal projective birational model $\pi\colon\cB'\to\cB_F$ such that the $b$ rational maps $\cB_F\dashrightarrow\bP^{n-1}$ defined by rows of $A$ extend to morphisms $\cB'\to\bP^{n-1}$. Denote by $\overline E_i$ the pullback of $\overline{\cO(1)}$ on $\bP^{n-1}$ to $\cB'$ via the morphism defined by $A_i$. Each $\overline E_i$ is nef and we have 
\begin{equation}\label{eq:row-bundle-height}
\widehat{\deg}(\overline H_F^{\,r}\overline E_i)/[F:K]\leqslant h(A_i)+\frac12\beta_S\log n.
\end{equation}

Put $H_A=\max_i h(A_i)$ and $q=(H_A+\beta_S\log n)/V$. For $p>0$, multiplication by the matrix defines a map between finite sets
\[ \widehat H^0(\mathcal B_F,p\overline H_F)^n\to\prod_{i=1}^b\widehat H^0(\mathcal B',p\overline H_F+\overline E_i(\frac12\log n)).\]
The twist $\frac12\log n$ accounts for the Cauchy--Schwarz inequality in each row. Each $\overline E_i(\frac12 \log n)$ is nef with smooth metric and has parameter at most $q$ in Lemma~\ref{lem:arithmetic-target}, by \eqref{eq:row-bundle-height}. The lemma therefore applies and we have, for $p\geqslant\max\{2,q\}$, 
\begin{align}
  \log\#(\text{source})&=n\Bigl(\frac{V_F}{D!}p^D+o_F(p^D)\Bigr),\label{eq:arithmetic-source-count}\\
  \log\#(\text{target})&\leqslant b\Bigl(\frac{V_F}{D!}(p+q)^D+C_Fp^{D-1}\log p\Bigr).\label{eq:arithmetic-target-count}
\end{align}
The first line is the Hilbert--Samuel formula for the fixed bundle $\overline H_F$ \cite[Corollary~5.5]{Moriwaki2009}. The second is uniform in all the rows and the blowup. In particular, neither error term depends on $A$.

Using $\log(1-t)\leqslant -t-t^2/2$ for $0<t<1$ and $\log(1+t)\leqslant t$ for $t\geqslant 0$, if $p\geqslant D\lambda q$ then
\[ \frac bn\Bigl(1+\frac qp\Bigr)^D\leqslant\Bigl(1-\frac1\lambda\Bigr)\Bigl(1+\frac1{D\lambda}\Bigr)^D\leqslant\exp\Bigl(-\frac1{2R^2}\Bigr)\eqqcolon 1-\eta<1.\]
Choose an integer $P=P(F,R)\geqslant 2$ so large that the error in \eqref{eq:arithmetic-source-count}, divided by $n$, has absolute value at most $\eta V_Fp^D/(4D!)$ for $p\geqslant P$, and $C_Fp^{D-1}\log p\leqslant\eta V_Fp^D/(4D!)$ there. Taking $p=P+\lceil D\lambda q\rceil$ makes the source strictly larger than the target. Thus there exist two distinct inputs with the same image. Their difference lies in $\widehat H^0(\mathcal B_F,p\overline H_F+\overline\cO(\log 2))^n$ and is killed by $A$. Trivializing $H_F^p$ at the generic point gives $0\neq z\in F^n$. Moreover, 
\[p\overline H_F+\overline\cO(\log 2)-\overline E_{z,\max}\]
is effective on some birational model $\cB''$ of $\cB_F$. Intersecting with $\overline H_F^{\,r}$ gives
\[\begin{split}
h(z)&=\widehat\deg(\overline H_F^{\,r}\overline E_{z,\max})/[F:K]\leqslant p\widehat\deg(\overline H_F^{\,r+1})/[F:K]+\beta_S\log 2\\
&= pV+\beta_S\log2 \leqslant D\lambda(H_A+\beta_S\log n)+(P+1)V+\beta_S\log2.
\end{split}\]
This proves \eqref{eq:arithmetic-siegel}. 
\end{proof}

The same result applies to a mixed polarization whose coordinate height satisfies
\[ c\,h_{\overline H}(w)\leqslant h(w)\leqslant C\,h_{\overline H}(w)\]
with fixed $0<c\leqslant C$, for every coordinate tuple $w$ over every finite extension. In that case one may take $\sigma_S=(C/c)D$. In particular, suppose the model line bundles $\overline H_1,\ldots,\overline H_r$ defining the polarization are all nef and big. By continuity of arithmetic volume \cite[\S2]{Ikoma2015}, there exist positive rational numbers $a_i,b_i$ such that $\overline H_i-a_i\overline H$ and $b_i\overline H-\overline H_i$ are big, hence pseudo-effective. Intersecting with $\overline E_{w,\max}$ associated to the tuple $w$ on some birational model $\cB'$ of some $\cB_F$ gives the comparison with $c=\prod_i a_i$ and $C=\prod_i b_i$, by \cite[Lemma~2.4(3)]{Ikoma2015}. 

\bibliographystyle{alpha}
\bibliography{references}

\end{document}